\documentclass[a4paper,twoside]{amsart}
\usepackage{amsmath}
\theoremstyle{plain}
\newtheorem{definition}{Definition}[section]

\newtheorem{lemma}[definition]{Lemma}
\newtheorem{theorem}[definition]{Theorem}

\newtheorem{corollary}[definition]{Corollary}
\newtheorem{proposition}[definition]{Proposition}

\newtheorem{notation}[definition]{Notation}
\newtheorem{remark}[definition]{Remark}
\theoremstyle{definition}

\newcommand\jcheck [1]{\check #1^j}

\theoremstyle{remark}

\begin{document}
\def\ZZ{Z\!\!\! Z}
\def\RR{I\!\! R}
\def\CC{I\!\!\! C}
\def\pf{\par\noindent {\it Proof: }}

\def\eop{\hfill {\rm QED.}\vskip 0.2cm}

\def\cht{\check\tau}
\def\chs{\check\sigma}
\def\chr{\check\rho}
\def\nn{\frac{a_1}{| X |}}
\def\n{\omega }
\def\p #1{p_{1,1}^{#1}}
\def\L{{\mathbf{A}}} 
\def\La{{\mathbf{L}}}
\def\k{{d}}
\def\gcud{\cdot \!\!>}
\def\lcud{\!\!\cdot}
\def\gcunor{\cdot \!\!\!>}
\def\lcunor{<\!\!\!\cdot}
\def\gcu{\ \cdot \!\!\!\!>}
\def\lcu{<\!\!\!\!\cdot\ }
\def\rk{{\rm rk}}
\def\qbin #1#2{\binom{#1}{#2}_{\!\! q}}
\def\bperp{{\perp\!\!\!\!\perp }}
\def\FA{\ \forall}
\def\direct{\mathrel\triangleright\joinrel\mathrel<}
\def\vaca{(q^{d-1}-1)}
\def\eldelta{{1_{_X}}}
\def\direct{\mathrel\triangleright\joinrel\mathrel<}

\def\newiote #1){\iota_{#1}}

%
%
\newenvironment{PalabrasClave}{
       \list{}{\advance\topsep by0.15cm\relax\small
       \leftmargin=1cm
       \labelwidth=0.35cm
       \listparindent=0.35cm
       \itemindent\listparindent
       \rightmargin\leftmargin}\item[\hskip\labelsep
                                     \bfseries Keywords:]}
     {\endlist}

\newenvironment{Numeros}{
       \list{}{\advance\topsep by0.15cm\relax\small
       \leftmargin=1cm
       \labelwidth=0.35cm
       \listparindent=0.35cm
       \itemindent\listparindent
       \rightmargin\leftmargin}\item[\hskip\labelsep
                                     \bfseries Mathematics Subject Classification:]}
     {\endlist}


\title[Lattices and Norton Algebras]{Lattices and Norton algebras of  \\
Johnson, Grassmann and Hamming graphs}

\thanks{Partially supported by Secyt-UNC, CIEM-CONICET, ANPCyT}

\author{C. Maldonado$^\dag$ and D. Penazzi*}
\address{\vbox{\tiny \hbox{*Facultad de Matem\'atica, Astronom\'\i a y F\'\i sica}
\hbox{Universidad Nacional de C\'ordoba,CIEM-CONICET}}\quad \vbox{\tiny \hbox{$^\dag$ Facultad de Ciencias Exactas F\'\i sicas y Naturales} \hbox{Universidad Nacional de C\'ordoba,CIEM-CONICET}}}

\begin{abstract}
To each of the Johnson, Grassmann and Hamming graphs we associate a lattice and characterize
the eigenspaces of the adjacency operator in terms of this lattice
. We also show that each level of the lattice induces in a natural way a tight frame for each eigenspace. For the most important eigenspace we compute explicitly the constant associated to the tight frame. Using the lattice we also give a formula for the product of the Norton algebra
attached to that eigenspace.
\end{abstract}
\maketitle
\begin{Numeros}
05E30,06D99,17D99
\end{Numeros}
\begin{PalabrasClave}
Johnson, Grassmann, Hamming, lattice, adjacency operator, tight frame, Norton algebra
\end{PalabrasClave}

\section{Introduction}\noindent

Distance regular graphs are important in Algebraic Combinatorics \cite{BCN}
and have been generalized into other 
combinatorial objects such as association schemes \cite{Godsil,MT}.
\!Some classical examples include the Johnson, Grassmann and Hamming graphs.
  Diverse algebras are associated to them, see for instance the Terwilliger Algebra in
 
  \cite{BST,C,CD,HKM,KLW,LMP2006,T}.    
Another algebra involved to such schemes is the Norton Algebra. In the 1970's Norton  constructed some commutative nonassociative algebras (called ``Norton Algebras" by Conway and by Smith in \cite{Smith}), whose automorphism groups contain finite groups generated by 3-transpositions, and in \cite{CGS} 
 this notion of algebra was applied
 to the case of an algebra  constructed on the eigenspaces of the adjacency operator of an association scheme.
(As is well known, related to this, Griess constructed the Monster simple group \cite{Griess} as the  automorphism group of a commutative nonassociative algebra of dimension 196883+1. This algebra is known as the Monster algebra but also as the Conway-Griess-Norton algebra.)

In a recent work \cite{LMP}, we have studied the Norton algebra (in the sense of \cite{CGS}) related to the dual polar graphs. While studying this problem we realized that the construccion of a lattice associated to these spaces was helpful and that it has some interesting properties of its own . In particular the eigenspaces of the adjacency operator of the graph can be reconstructed from the lattice (see below). We wanted to extend these lattice results to the case of the Johnson, Grassmann and Hamming graphs, since there are some technical differences between them and the dual polar graphs.
Using this framework, we also study their Norton Algebras.

Let $X$ be the set of vertices of these graphs.
The adjacency
operator $\L$ of the set of functions $\RR^{X}=\{ f:X
 \rightarrow \RR\}$ induced by the distance on the graph gives a  decomposition of $\RR^{X}$ into
eigenspaces of $\L$.

We construct a graded lattice associated to the graph and characterize the
eigenspaces of $\L$ in terms of this lattice (Theorem \ref{eigen}). 
 
 We show that the elements of each level of the lattice induces
 in a natural way a tight frame for each eigenspace (Theorem \ref{tightframe}). For references about the theory of finite normalized tight frames see for example \cite{BF, FMRHO, KC1, KC2, VW1, VW2}.

The eigenspace $V_1$ corresponding to the second largest eigenvalue of
$\L$ is of particular importance since one can reconstruct the
whole graph from the projections of the canonical basis onto it.
We explicitly compute the constant of the tight frame attached to
$V_1$.(Proposition \ref{lambda1})

We use these and other constants associated to the lattice to 
give a formula for the product of the Norton algebra
attached to $V_1$.(Theorem \ref{nortonproduct}). 

This article is organized as follows: In Section \ref{def} we give
definitions. In Section \ref{ldpg},  we define the lattice.
In section
\ref{embedding}, we give  a convenient description for the
eigenspaces $V_i$ of $\L$. The technical Proposition \ref{arribayabajo}
 is crucial for the proof of Theorem \ref{eigen}.
In Section \ref{Frames} we obtain tight frames and calculate the different \\ constants
associated to them for each of the cases Johnson, Grassmann and Hypercube.
In Section \ref{norton} we compute the Norton product using these constants.

\section{Definitions}\label{def}
\subsection{Distance regular graphs and their Adjacency algebras}\cite{BI}\label{drg}\noindent

Given  $\Gamma=(X,E)$ a graph with distance  $d(\ ,\ )$  we say
that it is distance regular if for any $(x,y) \in X\times X$ such
that $d (x,y)=h$ and for all $i,j\ge 0$ the cardinal of the set
$\{z \in X \mid
d (x,z)=i \ \mbox{and}\ d (y,z)=j \}$  is a constant denoted by
$p_{ij}^h$ which is independent of the pair $(x,y).$

Let  $\Gamma=(X,E)$ be a distance regular graph of diameter $d$.
Let $Mat_{X}(\RR)$ denote the $\RR$-algebra of matrices with real
entries, where the rows and columns are indexed by the elements of
$X$. For $0 \leq i \leq d$, the {\it ith adjacency matrix} of
$\Gamma$ is:
$ (A_i)_{xy} = \left \{ \begin{array}{ll}
1 & \mbox{if} \ d(x,y)=i \\
0 & \mbox{if} \ d(x,y)\neq i
\end{array}
\right. $.  It is easy to see that
the adjacency matrices satisfy:

 (i') $A_0=I$ where $I$ is the identity matrix;
(ii') $A_0+\dots +A_d=J$ where $J$ is the all $1' s$ matrix;
(iii') $ A_iA_j=\sum_{h=0}^d p_{ij}^h A_h \ (0\leq i,j \leq d)$;
(iv') ${A_i}^t=A_i$.
Thus $A_0,\dots, A_d$ form a basis for
a subalgebra $\mathcal{A}$ of $Mat_{X}(\RR)$ called the {\it adjacency algebra} of $\Gamma$.

Recall that there exists  a decomposition $\RR^X=\oplus_{j=0}^{d}
W_j$ where $\{W_j\}_{j=0}^{d}$ are common eigenspaces of
$\{A_i\}_{i=0}^{d}$. Let $p_i(j)$ the eigenvalue of $A_i$ on the
eigenspace $W_j$. By Proposition 1.1 of section 3.1 of Chapter III of \cite{BI},
 $\{A_i\}_{i=0}^{d}$ and the
eigenvalues $\{p_i(j)\}_{i,j=0}^{d}$ of a given $\Gamma$  satisfy: 
$ A_i=v_i(A_1) , \quad
p_i(j)=v_i(\theta_j),$ where $\theta_j=p_1(j)$, and
$\{v_i\}_{i=0}^{d}$ are polynomials of degree $ i$.
 We will order
the decomposition according to $\theta_0 > \theta_1
>...>\theta_d$.
In Theorem 5.1 of III.5 of \cite{BI}, one can find formulas for
the polynomials associated to each $\Gamma$.

\vspace{1em}

We will use the standar notations concerning the space of functions $\RR^X$:\\
i)\ $\mathbf{0}$ will denote the constant $\mathbf{0}(x)= 0 , \forall \ x \in
    X$,\ ii) the same for the constant $\mathbf{1}$, \ iii) $<f,g>:=\sum_{x\in
X}f(x)g(x)$,  \ iv)  $||f||^{2}:=<f,f>$, v) for $U\subseteq  \RR^X$,  $U^{\perp}=\{f \in \RR^X :\ <f,g>=0 \ \forall \ g \ \in U\}$. 

In addition, for ease of writing, we will use the following notation due to Iverson and Knuth (\cite{I,K}).
\begin{notation} (Iverson Bracket)\noindent

For any statement $P$, let $[P]=\begin{cases}  1 &
\mbox{if}\  P\ \mbox{is true.}\ \cr 0& \mbox{if}\ P \ \mbox{is
false.}\
\end{cases}
$

\end{notation}

\begin{definition}\label{definicionA}
Let $\L:\RR^X\mapsto \RR^X$ denote  the  adjacency operator
defined by $$\L(f)(x)=\sum_{y\in X}[d(x,y)=1]f(y)=\sum_{y\in X :\ d(x,y)=1}f(y)$$
\end{definition}

Observe that
$
\L (f)(x)
=\sum_{y\in X}(
A_1)_{xy}f(y)$. 
Then  $\{W_j\}_{j=0}^{d}$ are eigenspaces of $\L$ with $\theta_j$ as corresponding eigenvalues. 

$\L$ is symmetric and it holds that $<\L(f),g>=<f,\L(g)>$.

\subsection{Johnson, Grassman and Hypercube graphs}\label{dpg}\noindent

We define the distances regular graphs that we will use in the rest of the paper.\cite{BCN}

\vspace{.2em}

{\it Johnson graph}
The vertex set of  $J(n,k)=(X,E)$ ($2k \leq n $) is
the set of all $k$ -subsets of $[n]=\{1,...,n\}$, two vertices $x,y \in X $ being adjacent if and only if $|x\cap y|= k-1$ and as a consequence 
$d(x,y)=j \Leftrightarrow |x\cap y|= k-j. $  $J(n,k)$ has diameter $d=k$.

\vspace{.2em}

{\it Grassman graph} Let $V$ be an n-dimensional vector space over a field $F$ of $q$ elements. The vertex set of  $J_q(n,k)=(X,E)$ is 
the collection of linear subspaces of $V$ of dimension $k$. Two vertices $x,y \in X $ are adjacent if and only if 
$dim(x\cap y)= k-1$ and clearly
$d(x,y)=j \Leftrightarrow \ \dim (x\cap y)= k-j. $ $J_q(n,k)$ has diameter $d=k$.

\vspace{.2em}
{\it Hypercube graph} 
Take $S$ a set with two elements. The vertex set of \\ $H(n,2)=(X,E)$ is $S^n= \times _{i=1}^n S$ the cartesian product of $n$ copies of $S$, two vertices $x,y \in X $ being adjacent if and only if they differ precisely in one coordinate and therefore 
$d(x,y)=j \Leftrightarrow x  $ and $y$ differ precisely in $ j $ coordinates.
 $H(n,2)$ has diameter $d=n$.

\subsection{Lattice}\noindent

We recall the following definitions (see \cite{Stanley}) and in the next section we associate a lattice to each one of the distance regular graphs previously defined.

\begin{itemize}
\item A partial order is a binary relation "$\leq$" over a set P
which is reflexive, antisymmetric, and transitive.
\item A partially ordered set (poset) $(P,\leq)$ is a set $P$ with a partial order $\leq$.
\item A lattice $(\La,\leq,\wedge,\vee)$ is a poset  $(\La,\leq)$ in
which  every pair of elements $u, w \in \La$ has a least upper
bound and a greatest lower bound. The first is called  the join
and it is denoted by $u\vee w$ and the second is  called the meet
and it is denoted by $u\wedge w$.
\item A bounded lattice has a greatest (or maximum) and a least (or minimum) element, denoted $\hat{1}$ and $\hat{0}$ by convention.
\end{itemize}

\section{Lattice associated with Johnson, Grassman and Hypercubes graphs}\label{ldpg}

\subsection{}{\bf Johnson  graph} $J(n,k)$

\hspace{-.5cm} $\bullet$  For $j=0,...,k$ let
$\Omega_j$ be the
vertex set
 of $J(n,j)$ and $\Omega_{k+1}:=\{\hat 1\}$ where $\hat 1=[n]$
 
\hspace{-.5cm} $\bullet$ $\La=\cup_{\ell=0}^{k
+1}\Omega_\ell$ and for $x,y \in \La$ , \ $x \leq y \
\Leftrightarrow x \subseteq y $.

\hspace{-.5cm} $\bullet$  With that order $\La$ is a lattice with:    
 $x\wedge \ y=x\ \cap \ y,\ x\ \vee\ y=\left \{ \begin{array}{ll}
\!\! x \cup y &\mbox{if} \ |x \cup y |\leq k \\
\!\! \hat 1 &\mbox{if not}
\end{array}
\right. $

\subsection{} {\bf  Grassman graph $J_q (n , k)$} 

\hspace{-.5cm}$\bullet$  For $j=0,...,k$ let
$\Omega_j$ be the
vertex set
 of $J_q (n,j)$ and
$\Omega_{k+1}:=\left\{\hat 1\right\}$  where $\hat 1=V$

\hspace{-.5cm}$\bullet$ $\La=\cup_{\ell=0}^{k
+1}\Omega_\ell$ and for $x,y \in \La$ , \ $x \leq y \
\Leftrightarrow x \subseteq y $

\hspace{-.5cm}$\bullet$ With that order $\La$ is a lattice with:
    $$x\wedge y=x \cap y\ ,\ x\vee y=\left \{ \begin{array}{ll}
&  span \left\{x \cup y \right\}\ \mbox{if} \ dim(span\ \left\{x \cup y\right\} )\leq k \\
& \hat 1 \ \mbox{if not}
\end{array} 
\right. $$
\subsection{} {\bf The Hypercube $H(n,2)$} has as vertex set all words of length $n$ with symbols taken from
a set of $2$ elements. We will take as our set of two elements the
set $\{1,-1\}$, (instead of looking at words of $1$'s and $0$'s,
as is traditional).

Let $e_i$ be the vector with $n$ coordinates that has a $1$ in
position $i$ and $0$ elsewhere, and let $f_i=-e_i$. Then, each
word in $H(n,2)$ is simply a sum of some $e_i$'s ($i\in I$) and
some $f_j$'s,($j\in J$) with the only restrictions on $I $ and $J
$ that $I \cap J =\emptyset$ and $I \cup J =[n]$. 
Then the lattice associated is the following:

$\bullet$ For $0 \leq \ell \leq n $ we set $\Omega_\ell=\{\sum_{\{i \in I\}} e_i +\sum_{\{j \in J\}} f_j \ : \ I \cap J=\emptyset
    ,\ |I \cup J |= \ell\}.$    
    Given $ x \in \Omega_\ell,$ we represent $x=(I_x , J_x)$ where $x=\sum_{\{i \in
I_x\}} e_i +\sum_{\{j \in J_x\}} f_j $\\  $I_x \cap J_x=\emptyset$ and $ |I_x \cup J_x |= \ell$. 

For $x,y \in \cup_{\ell=0}^{n}\Omega_\ell ,$ we define  
 $x \leq y \ \Leftrightarrow \ I_x \subseteq J_x \
\mbox{and} \ I_y \subseteq J_y $. 

Observe that in the previous two cases 
$\Omega_\ell$ is  a member of the family of association schemes to which $X$
belongs even when $\ell < n$. This does not happen in this case since the words in $\Omega_\ell$
have $1,-1$ and $0$ in their entries while those in $X$ have only $1$'s and $-1$'s.

$\bullet$ Notice that  $\Omega_0:= \left\{(0,0,...,0)\right\}$ 
and we add a dummy element $\hat{1}$ above all other elements, defining $\Omega_{n+1}:=\left\{\hat 1\right\}$ 
, that is  $x\leq \hat 1 , \ \forall \ x  \in \cup_{\ell=0}^{n}\Omega_\ell$.

$\bullet$ With that order the set $\La=\cup_{\ell=0}^{n+1}\Omega_\ell$ is a lattice with:
$$
 x\wedge y=(I_x\cap I_y,J_x\cap J_y),\quad
x\vee y=\!\!\begin{cases} (I_x\cup I_y,J_x\cup J_y) & \mbox{if}\
(I_x\cup I_y)\cap (J_x\cup J_y)=\emptyset\cr \hat 1&
\mbox{otherwise}
\end{cases}$$
 (for $x,y \in \La-\left\{\hat 1\right\} $. Obviously $x\wedge\hat 1=x,\ x\vee\hat 1=\hat 1$)

\begin{definition}\noindent

 (1) Recall (\cite{Stanley}) that given 
 elements $u,w$ of a poset one says that 
  
 $u$ covers $w$; or $w$ is
covered by $u$, if $w < u$ but there is no $z$ such that $w<z<u$.
 
We denote it
by $u\gcunor w$ or $w\lcunor u$.

(2) A bounded lattice is ranked if the poset $\La$ is equipped with a rank function $ \rk:\La \rightarrow  \mathbb{Z}$ compatible with the ordering (so $\rk (u)\leq \rk(w)$ whenever $u \leq w$) and such that if $w$ covers $u$ then $ \rk(w)= \rk(u)+1$.
In our cases the lattices are clearly ranked, the rank of $w\in \Omega_\ell$ being $\ell$.

(3)An atom is an element that covers $\hat 0$ and a coatom is an element
covered by $\hat 1$. For $\Gamma=(X,E)$ any of the considered graphs,  the set
of atoms is $\Omega_1$ and the set of coatoms is $\Omega_d$ where ``$d$'' is the diameter of $\Gamma$. Since in fact the set of coatoms $\Omega_d$ is the set of vertices $X$ we will use both notations.

\end{definition}

Is not difficult to prove that for each graph defined above,
$(\La,\leq,\wedge,\vee)$ is a finite, bounded, ranked lattice with lowest element
$\hat 0$ and greatest element $\hat 1$.

\begin{remark}\label{loquenosolvidamos}
In all the cases:
$d(x,y)=j \Leftrightarrow x\wedge y\in \Omega_{d-j}.$ \ 
\end{remark}
\pf
\begin{eqnarray*}
d(x,y)=j\Leftrightarrow \ \left\{
\begin{array}{lll}
|x\cap y|=k-j &\mbox{for} \ J(n,k)\\
\dim (x\cap y)=k-j&\mbox{for} \ J_q(n,k)\\
x\mbox{ and } y\mbox{ have } n-j\\ 
\mbox{coordinates in common}&\mbox{for} \ H(n,2)\\
\end{array}\right\}\  
\Leftrightarrow x\wedge y\in \Omega_{d-j} \ .
\end{eqnarray*}
\eop

\begin{lemma}\label{prop}

The lattice $\La$  has the following properties:
\begin{enumerate}
    \item $\La$ is atomic, that is every element of the lattice is a join of atoms.
    \item $\forall u, w \in \La$ such that $   u\vee w\not=\hat 1\Rightarrow \rk(u)+\rk(w)=\rk(u\vee w)+\rk(u\wedge
    w)$ \label{modularity}
\end{enumerate}
\end{lemma}

\pf \noindent

(1)
In the Johnson  case, each element of the lattice is a set of
elements taken from $[n]$, so if $z=\{i_1,...,i_j\}$, then
$z=\{i_1\}\vee \{i_2\}\vee...\vee \{i_j\}$ is a join of atoms.

In the Grassman case, each element $z$ of the lattice is a
subspace of $GF(q)^n$, so taking a basis $\{v_1,...,v_j\}$ of $z$,
we obtain that $z=span(v_1)\vee span(v_2)\vee ...\vee span(v_j)$
is a join of atoms.

 In the Hamming case, $\hat 1= e_1 \vee f_1$ and  an element $z\neq \hat 1$ of the lattice is of the form
$z=\sum_{i\in I_z}e_i+\sum_{j\in J_z}f_j$ with $I_z\cap J_z=\emptyset$, so
$z=\bigvee_{i\in I_z}e_i\vee\bigvee_{j\in J_z}f_j$ is a join of
atoms.

(2)
In the Johnson case, $\rk(z)$  is the cardinality of $z$,
 so the
previous formula is simply the inclusion-exclusion formula for 2
sets: $|A|+|B|=|A\cup B|+|A\cap B|$.

In the Grassman, the rank of an element is
the dimension, so the formula is true because of the well known
identity $dim(V+W)=dim(V)+dim(W)-dim(V\cap W)$.

In the Hamming case, $\rk(z)=|I_z|+|J_z|$ 
so again the formula is true because of
the inclusion-exclusion formula for sets.
\eop
\begin{corollary}\label{first}\noindent

If $\tau $ and $\sigma $ are different atoms such that  
$\ \tau\vee\sigma\ne \hat 1 $, then $\rk(\tau\vee\sigma)=2 .$
\end{corollary}

\begin{lemma}\label{coverproperty}

Let $u$ and $w$ be elements of the lattice which are not coatoms, then
$(u\vee w \ \gcunor u , \ w)  \Rightarrow  (u ,  w \gcunor  u \wedge w)$.

 Reciprocally $(u ,  w\ \gcunor u\wedge
 w\ \mbox{and} \ u\vee w\ne\hat 1) \Rightarrow \ u\vee
 w \ \gcunor \ u , w$.
\end{lemma}

\pf In order to prove the first statement, observe that 
 $z \gcunor  w \Leftrightarrow z\ge w$   and $\ \rk(z)=\rk(w)+1$. 
 So, $u\vee w \gcunor u \mbox{ and} \ u\vee w \gcunor w \Leftrightarrow  \rk(u\vee w)=\rk(u)+1=\rk(w)+1$ (in
particular, we must have that $\rk(u)=\rk(w)$). Also, since $u$
and $w$ are not coatoms and $\rk(u\vee w)=\rk(u)+1$ we deduce that
$u\vee w\ne\hat 1$. Then, by  Lemma \ref{prop}
(\ref{modularity}),we get $\rk(u)+\rk(w)-\rk(u\wedge w)=\rk(u)+1$,
i.e., $\rk(w)=\rk(u\wedge w)+1$, which implies that $w \gcunor
u\wedge w$. The proof is similar for $u$.

 Reciprocally, if $u \gcunor u\wedge
 w$ and $w \gcunor u\wedge
 w \Rightarrow \rk(u)=\rk(w)=\rk(u\wedge w)+1$. Using  Lemma \ref{prop}
(\ref{modularity}) we get $\rk(w)=\rk(u)+\rk(w)-\rk(u\vee w)+1$
which implies $\rk(u)+1=\rk(u\vee w)$ and then that $u\vee w$
covers $u$ (and similarly $w$). \eop

\section{Description of the eigenspaces  using the associated lattice}\label{embedding}

\vspace{1em}

In this section let $\Gamma=(X,E)$ be any of the distance regular
graphs of diameter $d$ already defined ($J(n,k),J_q(n,k)$ or $H(n,2)$), together
with its associated decomposition:
$\ \RR^{X}=\oplus_{i=0}^{d} W_i , \ $ where $\{W_i\}_{i=0}^{d}$ are the
common eigenspaces of the adjacency matrices of $\Gamma$.

We will describe each of the eigenspaces $\{W_i\}_{i=0}^{d}$,
using the lattice previously defined. The description give us a recursive formulae for the eigenvalues $\left\{\theta_j\right\}_{j=0}^d$ associated to each graph $\Gamma=(X,E)$ defined in Subsection \ref{dpg}.

\begin{definition}\label{defal}
For $z\in\Omega_j$, define: 

$$a_j^{\ell}=\begin{cases} \{y\in \Omega_{\ell}:z\le y\} \
\mbox{if}\  j\le \ell  \cr \{y\in
\Omega_{\ell}:y\le z\}\ \mbox{if}\  j>\ell\cr \end{cases}\quad
a_j=\vert \{x\in X:z\le x\} \vert $$

Note that $a_j=a_j^{d}$ if $j\le d$, but $a_{d+1}=0$.

\end{definition}

The previous definitions seem to depend on $z$, we show next this is not so.

Recall that $\qbin {i}j$ is the number
of $j$-dimensional subspaces of $GF(q)^i$.
A formula is given by:
$\qbin {i}
j=\frac{[i]_q[i-1]_q ...[i-j+1]_q}{[j]_q[j-1]_q ...[1]_q}$
where 
$[i]_q=\left\{
\begin{array}{ccc} \frac{q^{i}-1}{q-1} & \forall \
i \geq 1\\ 0 & \forall \ i < 1 \end{array},\right .$ 

\begin{lemma}\label{al}
$$ 
a_j^{\ell}=\begin{cases}  {j \choose \ell} & \ \mbox{for}\ J(n,k)\cr \qbin {j}{\ell}&\ \mbox{for}\ J_q(n,k)\cr
 {j \choose \ell}  & \ \mbox{for}\ H(n,2) \cr
\end{cases}\ \mbox{if}\quad \ell \le j \
\quad
a_j^{\ell}=\begin{cases}  {n-j\choose l-j}& \ \mbox{for}\ J(n,k)\cr \qbin {n-j}{l-j}&\ \mbox{for}\ J_q(n,k)\cr
2^{\ell-j} {n-j \choose \ell-j} & \ \mbox{for}\ H(n,2) \cr
\end{cases}\ \mbox{if}\quad \ell \ge j \
$$
\end{lemma}

\pf
Given $z \in \Omega_j$ and $\ell \le j,$ we have to count the elements of the set\\
 $\{y\in \Omega_{\ell}:y\le z\}$. Looking at the construction of the lattice, the lemma follows straightforward in this case. 

In the case $j \le \ell$ the proof is also easy  for  $J(n,k)$ and $J_q(n,k)$. 
For the case $H(n,2)$, if we fix  $z\in\Omega_j$ and we count the elements of $\{y\in
\Omega_{\ell}:z\le y\}$, we have to choose $l-j$ coordinates from the $n-j$ not used by $z$, and we can fill each of them whith $1$ or $-1$.
\eop

\begin{definition}\label{defiota}
 \noindent

$\iota:\La  \rightarrow  \RR^X:z\mapsto \newiote z)$ is the map defined by $
\newiote z)(x)= [z\le x] \ \forall \ z \in \La ,\ x\in X .$

\end{definition}

\begin{lemma}\label{iotahat}\noindent

 For $\ j=0,1,...,d$ and  $ \forall \ z,y \in \La,$
 
\vspace{.5em}

$\begin{array}{llll} 

& i)\ \newiote z)=\mathbf{0} \Leftrightarrow z=\hat 1\quad &  ii)\ \newiote z)=\mathbf{1}\Leftrightarrow z=\hat 0 \

\quad iii)\ ||\newiote z)||^2=a_j\ \FA \ z\in \Omega_j
\\

\\
& iv)\
\newiote z)\newiote y)=\newiote z\vee y) 
& v)  <\newiote z),\newiote y)>=||\newiote z\vee y)||^2  \
\end{array}$
\end{lemma}
\pf \noindent

$i), ii)$ and $iii)$ are easy to prove. For $iv)$
$$
\newiote z)(x)\newiote y)(x)=[z\le x][y\le x]
=[z\vee y\le x]=\newiote z\vee y)(x)
$$

To prove $v)$,  observe that $$
<\newiote z),\newiote y)>=\sum_{x\in X}\newiote z)(x)\newiote y)(x)
=\sum_{x\in X}\newiote z\vee y)(x)
=\sum_{x\in X}(\newiote z\vee y)(x))^2=||\newiote z\vee  y)||^2 
$$
\eop
\begin{corollary}
\noindent

{\rm(i)} \ $z\vee y=\hat 1$ if and only if $\newiote z)$ and $\newiote y)$ are orthogonal to each other.

{\rm(ii)} If $z\vee y\in \Omega_j$, then $<\newiote z),\newiote y)>=a_j$

{\rm(iii)} If $\tau $ and $\sigma$ are both atoms then
$<\newiote \tau ),\newiote \sigma)>=\begin{cases}
a_1& {\rm if}\ \tau=\sigma \cr 0& {\rm if}\ \tau \vee
\sigma=\hat{1} \cr a_2& {\rm otherwise}
\end{cases}$

\end{corollary}

\subsection{A filtration for $\RR^{X}$.}

\begin{definition}\label{deflambdaj}\noindent

For $j=0,1,...,d$, let $\Lambda_j \subseteq \RR^{X}$ be the
subspace generated by $\{\newiote x)\}_{x \in \Omega_j}$.
\end{definition}
We want to show that $\Lambda_j\subseteq\Lambda_{j+1}$ . That is, they form a filtration for $\RR^{X}$. We need some tools first.

\begin{definition}\label{stars}\noindent
 Given $w\in \La$, let:
\begin{eqnarray*}
w^*=\sum_{v}[v\gcud w]\newiote v)\qquad
w_*=\sum_{v}[v\lcunor w]\newiote v)
\end{eqnarray*}
\end{definition}

\begin{lemma}\label{lemawstar}\noindent
Given $w\in  \Omega_{j} \subseteq \La$
\begin{eqnarray}
 w^* &=&c_{j} \ \newiote w) \ \mbox{where}\ c_{j} \ \mbox{only depends on}\  j=\rk(w) \\
 w_* &=&a_{j}^{j-1}\ \newiote w)+ \Phi_w  \ \mbox{where}\  \Phi_w:X\rightarrow \left\{0,1\right\}, \ \Phi_w(x)=\left[w\wedge x \in \Omega_{j-1}\right]
\end{eqnarray} 
\end{lemma}
\pf
(1) Given $x\in X$, We have:
 \begin{align*}
w^*(x)&=\sum_v[v\gcud w][v\le x]=\left|\{v\in \Omega_{j+1}:w\le v\le x\}\right|\\
&= \left\{
\begin{array} {cll} 0 & \mbox{if}\ w \not\le x \\
 \left\{
\begin{array} {cll} k-j & \mbox{for} \ J(n,k)\\
\left[k-j\right]_q & \mbox{for}\ J_q(n,k)\\
n-j & \mbox{for} \ H(n,2) \quad   \end{array}\right.
& \mbox{if}\ w \le x \\
\end{array}\right.
\end{align*} 
The last equality is easy to prove for Johnson and Grassman graphs. 
In the Hamming case, if  $w=\sum_{i\in I_w}e_i+\sum_{j\in
J_w}f_j$ and  $x=\sum_{i\in I_x}e_i+\sum_{j\in J_x}f_j$ to
build $v$ we need to add to the sum constituing $w$ one $e_i$ with
$i\in  I_x-I_w$ or else one $f_j$ with $j\in  J_x-J_w$. So: 
\begin{align*}
|\{v\in \Omega_{j+1}:w\le v\le x\}|&=| I_x-I_w|+| J_x-J_w|
\\
&=| I_x|+| J_x|-(|I_w|+|J_w|)
=n-j
\end{align*}
That is $w^* =c_{j} \ \newiote w)$ and the constant $c_{j}$ only depends on the $ \rk(w)=j$.

To prove the identity (2), given $x\in X$ we have 
\begin{align*}
 w_*(x)&=\sum_v[w\gcud v]\newiote v)(x)=\left|v : \ v\lcunor w \ \mbox{and} \ v\le x\right|\\
 &\stackrel{(\dag)}{=}\left\{ 
 \begin{array} {lll} 
& a_{j}^{j-1}& \mbox{if} \  w\le x \ (\mbox{equivalently} \  \rk(w\wedge x) = j) \\
& 1 & \mbox{if} \   \rk(w\wedge x) = j-1 \\
 & 0& \mbox{if}\ \rk(w\wedge x) <j-1\\
 \end{array} \right.
 \end{align*}

Thus  $w_* =a_{j}^{j-1}\ \newiote w)+ \Phi_w $.
The proof of  $(\dag)$ follows  from the fact in the case $ \rk(w\wedge x) < j-1$ there cannot be any such $v$. This is because if $v\le x$ and  $v\lcunor w$ then $v=v\wedge x \le w\wedge x$ and $ j-1=\rk(v)\le \rk(w\wedge x)$. 
\eop
\begin{corollary}\label{corolambdas}\ 
$\Lambda_0\subseteq\Lambda_1\subseteq...\subseteq \Lambda_{d}=\RR ^X$
\end{corollary}
\pf It follows from definition of  $ \Lambda_j$ and part (1) of previous lemma.\eop
\begin{definition}
Let   $V_0=\Lambda_0$ and $V_j=\Lambda_j\cap\Lambda_{j-1}^\perp
\quad j=1,...,d$.
\end{definition}
We have that $\Lambda_j=V_0\oplus V_1\oplus ...\oplus V_j = \Lambda_{j-1} \oplus V_j$.
We want to show that for $j=1,...,d$, $V_j\neq  \{0\}$, that is
$\Lambda_{j-1}  \neq  \Lambda_{j}$ . To prove this, we need more
lemmas.
Recall (Definition \ref{definicionA}) that the operator $\L$ is $\L(f)(x)=\sum_{y\in X }[d(x,y)=1]f(y)$.
\begin{lemma}
If $x\in X$, then $\L(\newiote x))=\sum_{y\in X }[d(x,y)=1]\newiote y)$.
\end{lemma}
\pf
Note that for $x,y\in X$ we have $\newiote x)(y)=[x=y]$ Thus, if $z\in X$:
\begin{align*}
\L(\newiote x))(z)&=\sum_{y\in X}[d(z,y)=1]\newiote x)(y)
=\sum_{y\in X}[d(z,y)=1][x=y]
=[d(z,x)=1]\\
&=\sum_{y\in X}[d(x,y)=1][y=z]
=(\sum_{y\in X}[d(x,y)=1]\newiote y))(z)
\end{align*}
\eop

\begin{lemma}\label{lemasumaiotas}
Let $x\in X$. Then\ 
$
x_*=\L(\newiote x))+\newiote x)a_d^{(d-1)}$.
\end{lemma}
\pf
Note that for $x,y \in X=\Omega_d$ we have by the previous proof
that\\
$\L(\newiote x))(y)=[d(x,y)=1]=[x\wedge y\in \Omega_{d-1}]=\Phi_x(y)$,
where $\Phi$ is as in Lemma \ref{lemawstar}.
Hence, the result follows from that lemma.
\eop
\begin{proposition}\label{arribayabajo}
For each $j<d$ there are constants $\alpha_j,\beta_j$ such that
if $w\in \Omega_j$ then:
$$
\sum_{u\gcunor w}u_*=\alpha_j \newiote w)+\beta_jw_*
$$
The constants are:

$\alpha_j=\left\{\begin{array}{ll}
n-2j& J(n,k)\\ \left[n-2j\right]_qq^j& J_q(n,k)\\ 2(n-j)& H(n,2) \end{array}\right.
,\ \beta_j=\left\{\begin{array}{ll}
c_{j-1}\quad  J(n,k) , J_q(n,k) \\ c_{j-1}-1 \quad H(n,2) 
\end{array}\right.$\\
where $c_j$ was given in the proof of Lemma \ref{lemawstar}.
\end{proposition}
\pf
\begin{eqnarray*}
\sum_{u\gcunor w}u_*&=&\sum_{u\gcunor w}\sum_{z\lcunor u}\newiote z)=\sum_{u\in \Omega_{j+1}}\sum_{z\in \Omega_j}[z\vee w\leq u]\ \newiote z)\\
&=&\sum_{u\gcunor w} \newiote w)+\sum_{u\in \Omega_{j+1}}\sum_{z\in \Omega_j-w}[z\vee w=u]\ \newiote z)\\
&{=}&a_j^{j+1}\newiote w)+\sum_{z\in \Omega_j}[\rk(z\vee w)=j+1]\newiote z)\quad (1)
\end{eqnarray*}
Similar arguments show that 
$\sum_{v\lcunor w}v^*=a_j^{j-1}\newiote w)+\sum_{z\in \Omega_j} [\rk(z\wedge w)=j-1]\newiote z)\quad (2)$

\vspace{.2cm}

For $J(n,k)$ and  $J_q(n,k)$ it follows from Lemma \ref{prop}  (\ref{modularity}) that $
 \rk(z\vee w)=j+1  \Leftrightarrow \rk( z\wedge w)=j-1
$
 so the sums at the rightmost side in (1) and (2) are equal and 
$$\sum_{u\gcunor w}u_*- \sum_{v\lcunor w}v^* =(a_j^{j+1}-a_j^{j-1})\newiote w)
$$
By Lemma \ref{lemawstar} $\sum_{v\lcunor w}v^*=\sum_{v\lcunor w}c_{j-1}\newiote v)=c_{j-1}w_*$, thus in the Johnson and Grassman cases we have
$\sum_{u\gcunor w}u_*=(a_j^{j+1}-a_j^{j-1})\newiote w)+c_{j-1}w_*$.
The values of $\alpha_j=a_j^{j+1}-a_j^{j-1}$ and $\beta_j=c_{j-1}$ follow from
Lemmas \ref{al} and \ref{lemawstar}.

The case of $H(n,2)$ is different because it is easy to prove that in this case 
\begin{eqnarray*}
&&\left\{z\in \Omega_j: \rk(z\wedge w)=j-1\right\}=\\
&& \left\{z\in \Omega_j \ : \ \rk(z\vee w)=j+1  \right\} \cup \left\{z\in \Omega_j: \rk(z\wedge w)=j-1\ \mbox{and} \ z\vee w=\hat{1}  \right\} ,
\end{eqnarray*}
\begin{eqnarray*} 
.^..\sum_{u\gcunor w}u_*- \sum_{v\lcunor w}v^* 
& =&(a_j^{j+1}-a_j^{j-1})\newiote w)\\
&&- \underbrace{\sum_{z\in \Omega_j}\left[ \rk(z\wedge w)=j-1\right] \left[ z\vee w=\hat{1}\right]\newiote z)}_{\Psi_w} .\ (3)
\end{eqnarray*}
Now given $x \in X$ we will evaluate $\Psi_w(x)$.
\begin{eqnarray*}
\Psi_w(x)
&=& \sum_{z\in \Omega_j}\left[ \rk(z\wedge w)=j-1\right] \left[ z\vee w=\hat{1}\right]\left[z \le x\right]\\
&=&\left| \left\{z\in \Omega_j: \rk(z\wedge w)=j-1 ,  \ z\vee w=\hat{1}  \ \mbox{and} \ z\le x\right\}\right|.
\end{eqnarray*}

Suppose that there is such a $z$. 
Such $z$ must be unique:
since $z\wedge w\in \Omega_{j-1}$ then it must be
$w=(z\wedge w)\vee \sigma$, where $\sigma\in \Omega_1$.
Similarly $z=(z\wedge w)\vee \tau$ for some $\tau$. Since $z\vee w=\hat 1$
it must be $\tau=-\sigma$. Thus $z$ is uniquely defined if it exists.

Moreover, 
since  $(z\le x \ \Rightarrow \  w \wedge z \le w \wedge x)$ then  we must have
$j-1 \leq \ \rk(w \wedge x) $. But  $\rk(w \wedge x) = j \Leftrightarrow w\leq x$ and since  $z\leq x$ this implies that $ z \vee w \leq x$ which is an absurd since $z \vee w=\hat{1}$. So $ \rk(w \wedge x) = j-1$ and $w \wedge x=w \wedge z$. Thus $w\wedge x \in \Omega_{j-1}$
must hold if $z$ exists.
Conversely, if $w\wedge x \in \Omega_{j-1}$ and $w=(x\wedge w)\vee \sigma$
then $z=(x\wedge w)\vee (-\sigma)$ satisfies all the conditions.
Thus, we conclude that $\Psi_w(x)=\left[w\wedge x \in \Omega_{j-1}\right]=\Phi_w(x)$.
Hence  by Lemma \ref{lemawstar} (2)  equation (3) becomes
\begin{eqnarray*}
\sum_{u\gcunor w}u_*- \sum_{v\lcunor w}v^* 
 &=&(a_j^{j+1}-a_j^{j-1})\newiote w)-\Phi_w \\
 &=&(a_j^{j+1}-a_j^{j-1})\newiote w)-\left(w_* - a_j^{j-1}\newiote w)\right)=a_j^{j+1} \  \newiote w) - w_*.
\end{eqnarray*}
As before, $\sum_{v\lcunor w}v^*=c_{j-1}w_*$, thus
$\sum_{u\gcunor w}u_*=a_j^{j+1} \  \newiote w)+(c_{j-1}-1) w_*$.
Again the values of $\alpha_j=a_j^{j+1}$ and $\beta_j=c_{j-1}-1$
follow from Lemmas \ref{al} and \ref{lemawstar}.
\eop
\begin{lemma}\label{mui}\noindent 

For  $0\leq j\leq d ,$  there exist a constant $
\lambda_j$ such that $\L(v)-\lambda_j v\in \Lambda_{j-1}$, $\forall\ v\in \Lambda_j.$
\end{lemma}
\pf It is enough to prove it for elements of the spanning
 set $ \left\{\newiote u)\right\}_{u\in \Omega_j}$.
The proof is by reverse induction on the levels of the lattice, starting at  $j=d$. 
The inductive hypothesis will be:

There exists constants $\lambda_{j}$ and $\nu_{j}$ 
 such that   
 $\L(\newiote u))=\lambda_{j}\newiote u)+
\nu_{j}u_*$ for all $ u\in \Omega_{j}. $ 
\vspace{.5em}

This will prove the lemma since by Definition \ref{stars}, $u_* \in \Lambda_{j-1}$.
 
The inductive hypothesis is true for $j=d$ since 
if $x\in \Omega_d$ then by  Lemma \ref{lemasumaiotas}
$\L(\newiote x))=-a_d^{(d-1)}\newiote x)+x_*$.
Now assume the hypotesis true for $j+1$ and 
let us prove it for $j$. Let $ w \in \Omega_j$.
 By Lemma \ref{lemawstar} $\newiote w)= \frac{1}{c_j} \
   w^* =\frac{1}{c_j} \sum_{u\gcud w}\newiote u)$, thus
\begin{align*} 
\L(\newiote w))&=\frac{1}{c_j} \sum_{u\gcud w}
\L(\newiote u)){=}\frac{1}{c_j} \sum_{u\gcud w} \left(
 \lambda_{j+1}\newiote u)+
\nu_{j+1}u_* \right) \\
&\stackrel{(\ref{arribayabajo})}{=}\frac{\lambda_{j+1}}{c_j} \ w^*+\frac{\nu_{j+1}}{c_j} \left(\alpha_j\newiote w)+\beta_jw_*\right)\\
&\stackrel{(*)}{=}\left(\lambda_{j+1}+\frac{\nu_{j+1}\alpha_j}{c_j}\right)\ \newiote w)+ \frac{\nu_{j+1}\beta_j}{c_j}w_*
\end{align*}
\eop
\begin{corollary}\label{Linv}
For $j=0,...,d$, $\Lambda_j$ are $\L$-invariant subspaces of
$\RR^{X}$.
\end{corollary}
\pf
 This follows directly by the previous Lemma and Corollary \ref{corolambdas}.
\eop

\begin{theorem}\label{eigen}\noindent
For $j =0,...,d$,  $V_j=\Lambda_j \cap \Lambda_{j-1}^{\perp}$ are the
eigenspaces $W_j$ of $\L$  given in \ref{drg}, in that order.
The corresponding eigenvalues $\theta_j$ 
are the $\lambda_j$'s
given in Lemma \ref{mui}.\end{theorem}
\pf
Take $v\in V_j\ (\subseteq \Lambda_j)$. By Lemma \ref{mui}
$\L(v)=\lambda_jv+v'$ with  $v'\in \Lambda_{j-1}$  and by Corollary
\ref{Linv} $\L(v')\in \Lambda_{j-1}$. Then by definition of $V_j$:
\begin{eqnarray*} 0&=&<v,\L(v')>=<\L(v),v'>=<\lambda_jv+v',v'>\\
&=&<\lambda_jv,v'>+<v',v'>=||v'||^2
\end{eqnarray*}
 thus  $\L(v)=\lambda_jv \ \forall \ v\in V_j$.
 Therefore, $\RR^{X}=\bigoplus_{j=0}^{d} V_j$ where each $V_j$ is either zero or an eigenspace of $\L$.
Since $X=\Omega_{d}$ is the set of vertices of  $\Gamma=(X,E)$;
a distance regular graph of diameter ${d}$;  there
are exactly $d+1$ eigenspaces of the adjacency matrix $A_1$,
therefore of $\L$. Thus each $V_j$ is indeed an eigenspace of $\L$
(hence $V_j\neq 0 \ \forall \ j$) and $\lambda_j$ are the eigenvalues
of $\L$. 

From the proof of Lemma \ref{mui} (identity $(*)$) we get that $\nu_d=1$ and the recursion $\nu_j=\frac{\nu_{j+1}\ \beta_j}{c_j} \ \forall j < d $. From the values of the constants, it follows  that $\nu_j=c_{j-1}$ in the Johnson and Grassman cases and $\nu_j=1$ in the Hamming case.
Therefore we conclude that $\lambda_j=\lambda_{j+1}+\alpha_j$ in the first two cases
and  $\lambda_j=\lambda_{j+1}+\frac{\alpha_j}{c_j}=\lambda_{j+1}+2$ in the latter case.
Therefore it is clear that $\lambda_0>\lambda_1>...>\lambda_d$. This imply (by
the ordering of 
\ref{drg}) that  
$\lambda_j=\theta_j$. 
\eop
\begin{remark}
From the recursion of the $\lambda$'s and the fact that $\lambda_j=\theta_j$, we obtain that the eigenvalues of $\L$ satisfy the following recursive formulae:
$$
\theta_d =\left\{\begin{array}{lll}
&- k&\mbox{for } J(n,k)\\ 
&- [k]_q& \mbox{for } J_q(n,k)\\
&- n&\mbox{for } H(n,2)
\end{array}\right.\qquad
\theta_j = \left\{\begin{array}{lll}
&\theta_{j+1}+n-2j & \mbox{for } J(n,k)\\
&\theta_{j+1}+[n-2j]_qq^j & \mbox{for } J_q(n,k)\\
&\theta_{j+1}+2& \mbox{for } H(n,2)
\end{array}\right.
$$
hence they are:
$$\theta_j=\left\{\begin{array}{lll}
&(k-j)(n-k-j)-j & \mbox{for } J(n,k)\\
&q^{j+1}[k-j]_q[n-k-j]_q-[j]_q & \mbox{for } J_q(n,k)\\
&n-2j & \mbox{for } H(n,2)
\end{array}\right.
$$
 The formulae for  $\theta_j$ can be founded in the literature (see Chapter 9 of \cite{BCN}). The proof above gives another way to compute them.
\end{remark}
\section{Tight Frames for the eigenspaces}\label{Frames}
In this section we will consider $\Gamma=(X,E)$ any of the graphs already defined,  $\La$ the
associated lattice described in Section \ref{ldpg} and \ 
$\RR^{X}=\oplus_{j=0}^{d} V_j$ \ the corresponding decomposition.
We will prove that each $\Omega_j$ induces a finite tight frame on each $V_j$ via the map defined in \ref{defiota}. We will give a formula for
the constants associated to these tight frames and in the case of  the eigenspace of the second largest
eigenvalue we will compute  explicitly the
constant associated. 

\begin{definition}\noindent

For $j=0,1,...,d$; let  $\pi_j$ be the orthogonal projection
 $\pi_j:\RR^{X}\rightarrow V_j$.
Then for each $u\in \Omega_j$, denote
 $\check u^j=\pi_j(\newiote u))$. Since the set $\left\{ \newiote u)\right\}_{u\in \Omega_j}$ span  $\Lambda_j$, the projections  $\left\{\check u^j\right\}_{u\in \Omega_j}$ span $V_j$. When it is obvious from the context we will denote it by $\check u$.
\end{definition}

\begin{proposition}\label{Uj}
For $j=0,1,...,d$, let  $U^{j} \in \RR^{X\times X}$ be the matrix
$$(U^{j})_{x,y}=(x,y)^j
=\sum_{u\in \Omega_j}\newiote u)(x)\newiote u)(y).$$
Then for every $j=0,...,d$, $\jcheck u$ is an eigenvector of $U^{j} $
with eigenvalue \\
$\mu_j=\sum_{i=0}^{d-j} a_{d-i}^j p_{i}(j)$
where $p_{i}(j)$ are the eigenvalues of $A_i$ (the $i$-th adjacency matrix of the graph) corresponding to the
eigenspace $V_j$.
\end{proposition}
 \pf
Let $(x,y) \in X\times X$ and $l=\rk(x\wedge y )$.
$$\begin{array}{rcll}
U^{j}_{x,y}&=&(x,y)^j
=\sum_{u\in \Omega_j}\newiote u)(x)\newiote u)(y)
=\sum_{u\in \Omega_j}[u\le x][u\le y]
=\sum_{u\in\Omega_j}[u\le x\wedge y]\\
&&\\
&=&\left\{\begin{array}{ccl}
 |\{u\in\Omega_j:u \le x\wedge y \}| \ &\mbox{if} &\ j\le l \\
 0 &\mbox{if}& \ j > l \end{array} \right. \qquad
= \left\{\begin{array}{ccl}
a_l^j \ &\mbox{if} &\ j\le l \\
0 &\mbox{if}& \ j > l \end{array} \right. \\
\end{array}$$

This and Remark \ref{loquenosolvidamos} shows that $U^{j}=\sum_{l=j}^{d}a_l^j  A_{d-l}.$
Then, as $\jcheck u \in V_j$ is an eigenvector of the
adjacency matrices, we have that for every $0\leq j \leq d$, $\jcheck u$ is an eigenvector of $U^{j} $
with eigenvalue $\mu_j=\sum_{l=j}^{d} a_l^j  p_{d-l}(j)$.

Making the change of variable $i=d-l$,
we have: $\mu_j=\sum_{i=0}^{d-j}a_{d-i}^j p_{i}(j)$.
\eop
\begin{definition}\noindent
Let $V$ be a finite vector space with inner product $<,>$. A finite tight frame on $V$ is a
finite set $F \subseteq V$ which satisfies the following
condition: there exists a non-zero constant $\mu$ such that:
$$ \sum_{v\in F} | <f,v>|^{2} = \mu \ \| \ f \|^{2}  \quad \forall \ f \in V.$$
\end{definition}
\begin{theorem}\label{tightframe}
For $j=0,...,d$ and 
for all $f\in V_j$, it holds 
$$\sum_{u\in\Omega_j}<\jcheck u,f>\jcheck u= \mu_j\ f$$ where $\mu_j$ is the eigenvalue of Proposition \ref{Uj}.
 
 In particular the set  $\{\jcheck u\}_{u\in\Omega_j}$ is a
finite tight frame for $V_j$. 

\end{theorem}
\pf
Let $w,v\in \Omega_j$.
\begin{align*}
<\mu_j \jcheck w,\jcheck v>&=\sum_{x\in X}\mu_j \jcheck w(x)\jcheck v(x)\stackrel{\ref{Uj}}{=}\sum_{x\in X}(\sum_{y\in X}(x,y)^j\jcheck w(y)\ )\jcheck v(x) \\
&=\sum_{x,y\in X}\sum_{u\in \Omega_j}\newiote u)(x)\newiote u)(y)\jcheck w(y)\jcheck v(x)\\
&=\sum_{u\in \Omega_j}<\newiote u),\jcheck w><\newiote u),\jcheck v>=\sum_{u\in \Omega_j}<\jcheck u,\jcheck w><\jcheck u,\jcheck v>\\ 
&=<\sum_{u\in \Omega_j}<\jcheck u,\jcheck w>\jcheck u,\jcheck v>.
\end{align*}
Since this is true for an arbitrary elements of the spanning set  $\left\{\check v^j\right\}_{v\in \Omega_j}$ of $V_j$ then $\mu_j \jcheck w=\sum_{u\in \Omega_j}<\jcheck u,\jcheck w>\jcheck u$, and again since this holds for arbitrary $w$ then it holds for any element of  $V_j$.
\eop
\subsection{Computation of $\mu_1$}
In the following we give a more explicit calculation of $\mu_1$, the constant associated to the tight frame corresponding to $V_1$; the second largest eigenspaces of $\Gamma=(X,E)$

\begin{proposition}\label{lambda1}
 $$
 \mu_1= \left\{\begin{array}{lll}
 {n-2 \choose k-1}  &\ \mbox{for} \ J(n,k) \\
 \qbin {n-2} {k-1} q^{k-1} &\ \mbox{for} \ J_q(n,k) \\
2^{n-1} & \ \mbox{for}  \ H(n,2) \end{array}\right.
$$
\end{proposition}
\pf
By Proposition \ref{Uj} we have that $\mu_1=\sum_{i=0}^{d-1} a_{d-i}^1 p_{i}(1)$.
One cand find the following formulae for $p_{i}(1)$ in pages 220 of \cite{BI} for $J(n,k)$; 262,263, 302 of \cite{BI} for $J_q(n,k)$ and 210 of \cite{BI} for $H(n,2)$;
$$
p_{i}(1)=
\left\{ \begin{array}{lll}&\sum_{t=0}^i (-1)^t 
{1 \choose t} {k-1 \choose i-t}{n-k-1 \choose i-t} & \mbox{for} \ J(n,k)\\
\\
&\left(1-\frac{\left[i\right]_q \left[n\right]_q}{\left[k\right]_q \left[n-k\right]_q \ q^i} \right)q^{i^2} \qbin k i \qbin {n-k} i &\mbox{for} \ J_q(n,k)\\
\\
& {n \choose i}- 2 {n-1 \choose i-1} &\mbox{for} \ H(n,2)
\end{array}\right.$$

Then by Proposition \ref{Uj}:
 
$\bullet$ For $ J(n,k)$: 
\begin{eqnarray*}
\mu_1&=&\sum_{i=0}^{k-1}a_{k-i}^1 \ p_{i}(1)=\sum_{i=0}^{k-1} (k-i) \left(\sum_{t=0}^i (-1)^t 
{1 \choose t} {k-1 \choose i-t} {n-k-1 \choose i-t} \right)\\
&=&\sum_{i=0}^{k-1} (k-i)\left({k-1 \choose i} {n-k-1 \choose i}-{k-1 \choose i-1} {n-k-1 \choose i-1}\right)
\end{eqnarray*}
Defining $b_i={k-1 \choose i} {n-k-1 \choose i}$\ (hence $b_i=0$ for $i < 0$).
\begin{eqnarray*}
&=&  \sum_{i=0}^{k-1} k\left(b_i-b_{i-1}\right)-\sum_{i=0}^{k-1} i\left(b_i-b_{i-1}\right) \\
&=&k\ b_{k-1}-\left(b_1-b_0+2(b_2-b_1)+3(b_3-b_2)+...\right)\\
&=&  k\ b_{k-1}-\left(-\sum_{i=0}^{k-1} b_i + k \ b_{k-1}\right)=\sum_{i=0}^{k-1}b_i={n-2 \choose k-1}
\end{eqnarray*}
(the last equation by  the Chu-Vandermonde identity: $ \sum_{i=0}^n {s \choose i}  {t \choose n-i} = {s+t \choose n}$).

$\bullet$ For $ J_q(n,k)$:
\begin{eqnarray*}
\mu_1&=&\sum_{i=0}^{k-1}a_{k-i}^1 \ p_{i}(1) \\
&=&\sum_{i=0}^{k-1} \left[k-i\right]_q \left(1-\frac{\left[i\right]_q \left[n\right]_q}{\left[k\right]_q \left[n-k\right]_q \ q^i} \right)q^{i^2} \qbin k i \qbin {n-k} i\\
&=&\sum_{i=0}^{k-1} \left[k-i\right]_q q^{i^2} \qbin k i \qbin {n-k} i - 
\sum_{i=0}^{k-1}  \frac{\left[k-i\right]_q \left[i\right]_q \left[n\right]_q q^{i^2} \qbin k i \qbin {n-k} i}{\left[k\right]_q \left[n-k\right]_q \ q^i} \\
\end{eqnarray*}
$$\mbox{Since:}\ \frac{\left[k-i\right]_q}{\left[k\right]_q}\ \qbin k i= \qbin {k-1} {i}=\qbin {k-1} {k-1-i},$$
$$\mbox{and:}\ \frac{\left[i\right]_q}{\left[n-k\right]_q} \qbin {n-k} i= \qbin {n-k-1} {i-1} = \qbin {n-k-1} {n-k-i}, $$
$$\mbox{then:}\
\mu_1
= \left[k\right]_q \sum_{i=0}^{k-1} \qbin {k-1} {k-1-i} \qbin {n-k} i \ q^{i^2}- \ 
\left[n\right]_q \sum_{i=0}^{k-1}   \qbin {k-1} i \qbin {n-k-1} {n-k-i}    q^{i^2-i}$$

Using $q$-Vandermonde: $\sum_{i}   \qbin m {k-i} \qbin {n} {i} q^{i(m-k+i)} = \qbin {m+n} k$  and  $[k]_q[n-1]_q=[n]_q[k-1]_q+[n-k]_qq^{k-1}$
we have: 
\begin{eqnarray*}
\mu_1&=&\left[k\right]_q \qbin {n-1} {k-1}- \ 
\left[n\right]_q \qbin {n-2} {n-k}=\left[k\right]_q \qbin {n-1} {k-1}-\left[n\right]_q\qbin {n-2} {k-2}  \\
&=& \qbin {n-2} {k-2}\left(\left[k\right]_q \frac{\left[n-1\right]_q}{\left[k-1\right]_q}-\left[n\right]_q\right)\\
&=& \qbin {n-2} {k-2}\ \frac{\left[n-k\right]_q}{\left[k-1\right]_q}\ \ q^{k-1}
=\qbin {n-2} {k-1}\  \ q^{k-1}\\
\end{eqnarray*}

$\bullet$ For $ H(n,2)$,
\begin{eqnarray*}
\mu_1&=&\sum_{i=0}^{n-1}a_{n-i}^1 \ p_{i}(1) =\sum_{i=0}^{n-1} (n-i) \left( {n \choose i}- 2 {n-1 \choose i-1} \right)\\
&=&\sum_{i=0}^{n-1} (n-i) \ {n \choose i}- 2 \sum_{i=0}^{n-1} (n-i) {n-1 \choose i-1} \\
 &=&\sum_{i=0}^{n-1} n \ {n-1 \choose i}- 2 \sum_{i=1}^{n-1} (n-1) {n-2 \choose i-1} 
 \\
 &=& n 2^{n-1}-2(n-1)2^{n-2}= 2^{n-1}
 \end{eqnarray*}
\eop
\section{Application: Norton product on $V_1$}\label{norton}
Given the decomposition $\RR^{X}=V_0\oplus V_1\oplus...\oplus
V_d$, in this section we describe the product of a Norton algebra
attached to the eigenspace $V_1$.

\begin{definition}
The  Norton algebra on $V_1$ is the algebra given by the product
$f\star g=\pi_1(fg)$ for $f,g\in V_1$.

\end{definition}

It is easy to check that it is a commutative, nonassociative algebra. We want to compute the $\star $ product in $V_1$ for the graphs concerning on this paper. Since
$\Lambda_1=span\{\newiote \tau):\tau\in \Omega_1 \}$ the set
$\{\check \tau\}_{\tau\in \Omega_1}$ spans $V_1$ and we have proved in Theorem \ref{tightframe} that they are a a finite tight frame for $V_1$. We will describe $\cht\star \chs$ in a simplified form using such a frame.

For this we need the following results.

\begin{lemma}\label{taucheck} 
For all $\tau\in\Omega_1$, \ $\cht=\newiote \tau)- \frac{a_1}{| X |}
\mathbf{1} $ with $a_1$  given in Lemma \ref{al}.
\end{lemma}

\pf Recall that $<\newiote \hat{0})>=\Lambda_0 \ \subseteq
\Lambda_1=<\{\newiote \tau)\}_{\tau \in \Omega_1}>$, and
$\Lambda_1=\Lambda_0\oplus V_1$. Since $\forall \ \tau\in\Omega_1 ,\ \cht=\pi_1(\newiote \tau))
\in V_1$, we
have $\cht= \newiote \tau) -t.\mathbf{1}$ for some $t \in \RR$.\\
From the fact that $<\cht,\mathbf{1}>=0$ we conclude 
$t=\frac{<
 \newiote \tau) ,\mathbf{1}>}{||\mathbf{1}||^2}=\frac{\sum_{x\in
X}[\tau\leq x]}{ |X|}=\nn$.
\eop
\begin{proposition}\label{pitheorem}
Let $h \in \RR^{X}$, then \ \
$\pi_1(h)=\sum_{\tau\in \Omega_1}\frac{< \newiote \tau) ,h>}{\mu_1}
\check \tau .$
\end{proposition}

\pf 
\begin{align*}
\mu_1\pi_1(h)&\stackrel{\ref{tightframe}}{=}\sum_{\tau\in \Omega_1}< \check \tau ,\pi_1(h)>\check \tau \stackrel{\dag}{=}\sum_{\tau\in \Omega_1}< \check \tau ,h>\check \tau 
\stackrel{\ref{taucheck}}{=}\sum_{\tau\in \Omega_1}< \newiote \tau)-\frac{a_1}{| X |} \mathbf{1}
,h>\check \tau\\
&=\sum_{\tau\in \Omega_1}< \newiote \tau), h>\check
\tau -<\frac{a_1}{| X |} \mathbf{1} ,h> \sum_{\tau\in
\Omega_1}\check \tau
\stackrel{*}{=}\sum_{\tau\in \Omega_1}< \newiote \tau), h>\check
\tau
\end{align*}
($\dag$ holds since $h=\pi_0(h)+\pi_1(h)+...+\pi_d(h)$, \ $\pi_i(h)\in V_i $ and 
$<V_i,V_j>=0 \ \forall \ i\neq j$; 
$*$ holds since $\sum_{\tau\in \Omega_1} \newiote \tau) \in \Lambda_0\Rightarrow \sum_{\tau\in
\Omega_1}\check \tau=\mathbf{0}$)
\eop
\begin{lemma}\label{taustarsigma}
$$\cht\star\chs=
\pi_1\left( \newiote \tau \vee \sigma)\right)-\nn (\cht+\chs)
$$
\end{lemma}
\pf Recall  by Lemma
\ref{taucheck} \ $\cht= \newiote \tau) -\frac{a_1}{| X |}\mathbf{1}$. Then: 
\begin{eqnarray*}
\cht \star \chs
&=&\pi_1 (\cht\  \chs)= \pi_1\left( ( \newiote \tau) -\frac{a_1}{| X |}\mathbf{1})(\newiote \sigma)-\frac{a_1}{| X |}\mathbf{1} ) \right)\\
&=&\pi_1\left(  \newiote \tau)  \newiote \sigma)-\frac{a_1}{| X |}\left(  \newiote \tau)  +\newiote \sigma) \right)+(\frac{a_1}{| X |})^2\ \mathbf{1}  \right)\\
&=&\pi_1\left(\newiote \tau)  \newiote \sigma)\right)-\frac{a_1}{| X
|}\pi_1\left(\newiote \tau)  +\newiote \sigma) \right)+(\frac{a_1}{| X
|})^2\pi_1(\mathbf{1})\\
&=&\pi_1\left( \newiote \tau \vee \sigma)\right)-\nn (\cht+\chs)\quad
\end{eqnarray*}
\eop
\begin{lemma}\label{rhotausigma}
If $a_j$ are as in Lemma
\ref{al}, then $<\newiote \rho), \newiote \tau \vee \sigma)
>=a_{\rk(\rho\vee\tau\vee\sigma)}.$
\end{lemma}
\pf
\begin{eqnarray*}
<\newiote \rho), \newiote \tau \vee \sigma) >&=&
\sum_{x\in X}\newiote \rho)(x) \newiote \tau \vee \sigma) (x)=\sum_{x\in X}[\rho\leq x][\tau \vee \sigma\leq x]\\
&=&\sum_{x\in X}[\rho\vee\tau \vee\sigma\leq x]=|\{x\in X:\rho\vee\tau \vee\sigma\leq x\}|=a_{\rk(\rho\vee\tau\vee\sigma)}\\
\end{eqnarray*}\eop
\def\sumarara{\sum_{\rho \le \tau\vee \sigma} \chr }
\begin{theorem}\label{nortonproduct}\noindent
{For} $H(n,2)$, $\cht\star\chs=\mathbf{0}$. 

$\cht\star\cht=(1-\frac{2k}{n})\cht$ in the Johnson case and $(1-\frac{2\left[k\right]_q}{\left[n\right]_q})\cht$ in the Grassmann case,
while for $\tau\neq\sigma$:
$$
\cht\star\chs=\begin{cases} 
\frac{2k-n}{n \ (n-2)} (\cht+\chs)
&\mbox{For} \ J(n,k)\cr
\cr
  - \frac{[k]_q}{[n]_q}(\cht+\chs)+
\frac{[k-1]_q}{q[n-2]_q}\sumarara &
\mbox{For} \ J_q(n,k) \cr
\end{cases} 
$$
\end{theorem}
\pf \noindent
By Lemma \ref{taustarsigma},
if $\tau=\sigma $ we have that $\cht\star\cht= \cht-2 \nn \cht
$. Replacing $a_1$ by Lemma \ref{al} the formulae follow straighforward for all the  graphs.

For the case $\tau\neq\sigma $ we will use the notation
$\Psi_j=\{\rho\in \Omega_1: \rk(\rho\vee\tau\vee\sigma)=j\}$.
Using 
 Lemmas \ref{taustarsigma} and \ref{rhotausigma}:
\begin{eqnarray*}
\cht\star\chs&=& -\nn
(\cht+\chs)+\pi_1(\newiote \tau \vee \sigma))\\
&=& -\nn
(\cht+\chs)+
\sum_{\rho\in\Omega_1}\frac{<\newiote \rho), \newiote \tau \vee \sigma) >}{\mu_1}\chr = -\nn
(\cht+\chs)+\sum_{\rho\in\Omega_1}\frac{a_{\rk(\rho\vee\tau\vee\sigma)} }{\mu_1} \chr\\
 &=& -\nn
(\cht+\chs)+\frac{a_2}{\mu_1}\sum_{\rho\in\Psi_2}\chr+\frac{a_3}{\mu_1}\sum_{\rho\in\Psi_3}\chr+\mathbf{0}
\end{eqnarray*}
 The last zero since  $a_{d+1}=0$. Also, since:
 $$\sum_{\rho\in\Psi_3}  \chr=\sum_{\rho\in\Omega_1}  \chr-\sum_{\rho\in\Psi_2}  \chr-\sum_{\rho\in\Psi_{d+1}}  \chr$$
and 
$\sum_{\rho\in\Omega_1} \newiote \rho) \in \Lambda_0\Rightarrow \sum_{\rho\in\Omega_1} \chr=\mathbf{0} $
we have then:
$$(\diamondsuit)\quad \cht\star\chs= -\nn
(\cht+\chs)+ \frac{a_2-a_3}{\mu_1}\sum_{\rho\in\Psi_2}\chr-\frac{a_3}{\mu_1}\sum_{\rho\in\Psi_{d+1}}  \chr$$

\noindent Then, we have, in each case:

$\bullet$ For $J(n,k)$, $\Psi_2=\{\tau,\sigma\}$ and $\Psi_{d+1}=\emptyset$. 
In this case then $(\diamondsuit)$ becames:
\begin{eqnarray*} 
\cht\star\chs &=&
 -\nn
(\cht+\chs)+\frac{a_2-a_3}{\mu_1}(\cht+\chs)\\
 &=&\left(-\frac{{n-1 \choose k-1}}{{n \choose k}}+\frac{ {n-2 \choose k-2}-{n-3 \choose k-3}}{{n-2 \choose k-1}}\right)(\cht+\chs)\\
  &=& \left(-\frac{k}{n}+\frac{k-1}{n-2}\right) (\cht+\chs)=
 \frac{2k-n}{n\ (n-2)} (\cht+\chs)
\end{eqnarray*}
$\bullet$ For $J_q(n,k)$, $\Psi_2=\{\rho\in\Omega_1:\rho \le \tau\vee \sigma\}$ and $\Psi_{d+1}=\emptyset$. 

Then $(\diamondsuit)$ becames:
$\cht\star\chs =
 -\nn
(\cht+\chs)+\frac{a_2-a_3}{\mu_1} \sum_{\rho \le \tau\vee \sigma} \chr 
$
\def\sumarara{\sum_{\rho \le \tau\vee \sigma} \chr }

Recall that in this case, $a_j=\qbin {n-j}{k-j}$ (Lemma \ref{al}), $|X|=\qbin nk$ and $\mu_1=\qbin {n-2} {k-1} q^{k-1}$ (Proposition \ref{lambda1}). 
\begin{eqnarray*}
\mbox{Thus}\quad \cht\star\chs &=&- \frac{[k]_q}{[n]_q}(\cht+\chs)+
\frac{1-\frac{[k-2]_q}{[n-2]_q}}
{\frac{[n-k]_q}{[k-1]_q}\ q^{k-1}}\sumarara\\
&=&- \frac{[k]_q}{[n]_q}(\cht+\chs)+
\frac{([n-2]_q-[k-2]_q)[k-1]_q}{[n-2]_q[n-k]_qq^{k-1}}\sumarara\\
&=&- \frac{[k]_q}{[n]_q}(\cht+\chs)+
\frac{\left((q^{n-2}-1)-(q^{k-2}-1)\right)(q^{k-1}-1)}{(q^{n-2}-1)(q^{n-k}-1)q^{k-1}}\sumarara\\
&=&- \frac{[k]_q}{[n]_q}(\cht+\chs)+
\frac{q^{k-2}\left(q^{n-k}-1)\right)(q^{k-1}-1)}{(q^{n-2}-1)(q^{n-k}-1)q^{k-1}}\sumarara\\
&=&- \frac{[k]_q}{[n]_q}(\cht+\chs)+
\frac{[k-1]_q}{[n-2]_q\ q}\sumarara
\end{eqnarray*}

$\bullet$ For $H(n,2)$  $\Psi_2=\{\tau , \sigma\}$ and $\Psi_{d+1}=\{-\tau , -\sigma\}$ and it holds that
$\forall \rho \in  \Omega_1$ $\quad \newiote \rho)+\newiote -\rho)= \mathbf{1} $ therefore: $ \check{\rho}+\check{(-\rho)}=\mathbf{0}\ \forall \rho \in  \Omega_1$,\  i.e. $\check{(-\rho)}=-\check{\rho}$. Then:
\begin{eqnarray*}
\cht\star\chs &=&
 -\nn
(\cht+\chs)+\frac{a_2-a_3}{\mu_1} (\cht+\chs)-\frac{a_3}{\mu_1} \left(\check{(-\tau)}+\check{(-\sigma)}\right)\\
&=&-\frac{2^{n-1}}{2^{n}}(\cht+\chs)+\frac{2^{n-2}-2^{n-3}}{2^{n-1}}(\cht+\chs)-\frac{2^{n-3}}{2^{n-1}}\left(-\cht-\chs\right)\\
&=&-\frac{1}{2}(\cht+\chs)+\frac{1}{2}\left(\cht+\chs\right)
=\mathbf{0}
\end{eqnarray*}
\eop
\begin{remark}
The fact that in the Hamming case 
the Norton product reduces to zero can also be deduced from
Theorem 5.2 of \cite{CGS} since it can be shown that the ``Krein parameters" $q_{1,1}^1$ are $0$ in this case.
\end{remark}

\section{Conclusion}
For each of the Johnson, Grassmann and Hamming graphs we constructed a ranked (finite) lattice which we embeed into $\RR^X$ (Definition \ref{defiota}). For the levels $\Omega_j$ the corresponding  embeedings $\Lambda_j$ in $\RR^X$ are shown to be a filtration, and  we characterized
the eigenspaces $W_j$ of the adjacency operator in terms of these $\Lambda_j$s.
(Theorem \ref{eigen}).
 We also show that each $\Omega_j$  induces in a natural way a tight frame for each eigenspace. Using the lattice we  give a formula for the product of the Norton algebra
attached to $W_1$.

\end{document}